\documentclass[11pt]{amsart}

\usepackage{amsmath,amsthm,amssymb,amsfonts,mathtools}
\usepackage{enumitem}
\usepackage{xcolor}
\usepackage{tikz}
\usetikzlibrary{arrows.meta,calc,positioning,patterns,backgrounds}
\usepackage[colorlinks=true,linkcolor=blue,citecolor=green,
  urlcolor=blue]{hyperref}
\hypersetup{
  pdftitle={Twisted Diophantine approximation: asymptotic and uniform results},
  pdfauthor={Taehyeong Kim and Vasiliy Neckrasov}
}

\setenumerate{leftmargin=*}
\setitemize{leftmargin=*}
\allowdisplaybreaks
\numberwithin{equation}{section}
\makeatletter
\renewcommand{\thmhead@plain}[3]{%
  \thmname{#1}\thmnumber{\@ifnotempty{#1}{ }#2}%
  \thmnote{ {\the\thm@notefont #3}}%
}
\let\thmhead\thmhead@plain
\makeatother
\newtheorem{thm}{Theorem}[section]
\newcommand{\R}{\mathbb R}
\newcommand{\Z}{\mathbb Z}
\newcommand{\N}{\mathbb N}
\newcommand{\T}{\mathbb T}
\newcommand{\Mat}{\mathrm{Mat}}
\newcommand{\bq}{\mathbf q}
\newcommand{\bp}{\mathbf p}
\newcommand{\bb}{\mathbf b}
\newcommand{\bx}{\mathbf x}

\newcommand{\Sing}{\mathrm{Sing}}
\newcommand{\Emass}{\operatorname{Emass}}
\newcommand{\e}{\mathrm e}
\renewcommand{\le}{\leqslant}
\renewcommand{\ge}{\geqslant}

\colorlet{figink}{black}
\colorlet{figmuted}{black!60}
\colorlet{figline}{black!35}
\colorlet{figorbit}{blue!65!black}
\colorlet{figwindow}{orange!80!black}

\title[Twisted Diophantine approximation]
{Twisted Diophantine approximation:\ asymptotic and uniform results}
\author{Taehyeong Kim}
\address{Department of Mathematics, Brandeis University,
Waltham, MA 02453, USA}
\email{taehyeongkim@brandeis.edu}
\author{Vasiliy Neckrasov}
\address{Department of Mathematics, Brandeis University,
Waltham, MA 02453, USA}
\email{vneckrasov@brandeis.edu}
\date{}

\begin{document}

\begin{abstract}
{Fix an $m\times n$ matrix $A$; twisted Diophantine
approximation studies how well we can approximate the affine form
$A\bq-\bb$ by an integer $m$-vector, depending on the choice of
$m$-vector $\bb$.}
We report on two companion papers in twisted Diophantine approximation{,
covering two setups: asymptotic and uniform.}
{The} asymptotic paper gives an exact
{Khintchine-type} Lebesgue zero--one law,
a null--conull criterion for badly approximable shifts, and
Hausdorff-measure tests yielding exact power-law dimensions under a
critical divergence condition. The uniform paper computes
{Hausdorff} dimensions {of sets of uniformly
approximable vectors.} It also gives
Hausdorff zero--full laws, uniform-exponent spectra, and a sharp escape-of-mass bound. In both cases the results are established for arbitrary approximation functions from suitable regularity classes, and the arguments use a homogeneous lattice count involving all successive minima to record clustering in every independent direction.
\end{abstract}

\maketitle

\section{Introduction}
\label{sec:introduction}

\subsection{Diophantine approximation: asymptotic and uniform}

Let $m,n$ be positive integers, set $d=m+n$, and write
$\Mat_{m,n}$ for the space of real $m\times n$ matrices.
{Classical theory of homogeneous Diophantine approximation asks,
for $A\in\Mat_{m,n}$ and a positive non-increasing function
$\psi:\R_+\to\R_+$, where $\R_+=(0,\infty)$,
whether there are integer vectors $\bp\in\Z^m$ and $\bq\in\Z^n$ satisfying}
\begin{equation}\label{eq:dirichlet}
 |A\bq-\bp|\le\psi(T),\qquad 0<|\bq|\le T
\end{equation}
at scale $T>0$. Here $T$ bounds the size of $\bq$, while
$\psi(T)$ specifies the allowed approximation error; throughout,
$|\cdot|$ denotes the supremum norm.

Dirichlet's theorem (1842) states that \eqref{eq:dirichlet} is solvable
for every $A\in\Mat_{m,n}$ and every $T>1$ when $\psi(T)=T^{-n/m}$.
{For a general $\psi$, we may ask two types of questions:
\emph{do solutions exist for an unbounded set of $T$?} Or
\emph{do solutions exist for all sufficiently large $T$?}} In the former case
we call $A$ \emph{$\psi$-approximable}, and in the latter
\emph{$\psi$-Dirichlet}\footnote{Sometimes the definitions of these sets
use a different normalization; see \cite{MN25,KW19}.}; these are the asymptotic and uniform
conditions, respectively \cite[Section~1.1]{KW19}.

The metric theory asks how large these two classes of matrices
are. For the asymptotic condition, the Khintchine--Groshev theorem
\cite{Khi26,Gro38} states that the set of $\psi$-approximable
matrices is null or conull according as
\begin{equation}\label{eq:classical-sum}
 \sum_{q=1}^\infty q^{n-1}\psi(q)^m,
 \quad\text{equivalently,}\quad
 \sum_{k=0}^\infty 2^{kn}\psi(2^k)^m
\end{equation}
converges or diverges; see \cite[Section~1]{BV10} for this
formulation. For the uniform condition, the set of $\psi$-Dirichlet
matrices is likewise null or conull under suitable regularity
assumptions on $\psi$. This zero--one law was established by
Kleinbock and Wadleigh for $m=n=1$ \cite{KW18}, while in
arbitrary dimensions the partial result of Kleinbock,
Str{\"o}mbergsson and Yu \cite{KSY22} was completed by
Str{\"o}mbergsson and Yu \cite{SY26}.

\subsection{Inhomogeneous approximation}

\begingroup
\emergencystretch=1em
{In the setup of inhomogeneous Diophantine approximation,
we replace the linear form $A\bq$ by an affine form $A\bq-\bb$,
where $\bb\in\R^m$. We study the system}
\begin{equation}\label{eq:system}
 |A\bq-\bp-\bb|\le\psi(T),\qquad |\bq|\le T.
\end{equation}
\endgroup
The same distinction between asymptotic and uniform approximation
now applies to pairs: we call $(A,\bb)$ \emph{$\psi$-approximable} if
\eqref{eq:system} is solvable for unbounded $T$, and
\emph{$\psi$-Dirichlet} if it is solvable for all
sufficiently large $T$, with $\bq=\mathbf0$ allowed here.
In the inhomogeneous setting, Dirichlet's bound is no longer
valid for every pair $(A,\bb)$.
Indeed, if $A$ has integer entries and $\bb\notin\Z^m$, then
\eqref{eq:system} has no solution once $\psi(T)$ is small enough.

The metric results depend on which parameters are allowed to vary.
When a nonintegral shift $\bb$ is fixed and $A$ varies, the
asymptotic problem retains the classical measure law: the set of
$\psi$-approximable matrices is null or conull according as the
series in \eqref{eq:classical-sum} converges or diverges
\cite{Szu58,Sch64,Spr79}. With both $A$ and $\bb$ varying, the
set of $\psi$-Dirichlet pairs is also null or conull by the
zero--one law of Kleinbock and Wadleigh \cite{KW19}.
The first-named author and Kim further studied the complementary non-Dirichlet sets,
establishing Hausdorff-measure laws for pairs and for fixed
$\bb\notin\Z^m$ when $\psi$ is decreasing and tends to zero
\cite{KK22}. In these results, the matrix $A$ remains a variable.
The question addressed in this report is how the two approximation
conditions behave when $A$ is fixed and only the shift $\bb$ varies.

\subsection{{Twisted approximation}}

This choice is usually called \emph{twisted approximation}.
From now on, we assume $\psi(T)\to0$ as $T\to\infty$.
Since integer translates of $\bb$
{do not affect solvability of \eqref{eq:system}}, shifts naturally
belong to $\T^m=\R^m/\Z^m$; we write
$\langle\bx\rangle=\bx+\Z^m$ for the image of $\bx\in\R^m$.
For a fixed $A$, let $W_A(\psi)$ and $D_A(\psi)$ be the sets of
shifts for which \eqref{eq:system} is solvable for unbounded $T$
and for all sufficiently large $T$, respectively.\footnote{{For
non-increasing $\psi$ tending to zero, this convention differs from that in
\cite{KNA26} only on the countable orbit
$\{\langle A\bq\rangle:\bq\in\Z^n\}$ and boundaries of balls
$B(\langle A\bq\rangle,\psi(|\bq|))$ with $\bq\ne\mathbf0$.
These sets are Lebesgue null.}}

\pagebreak[3]
{We briefly recall the main metric results for the asymptotic
sets $W_A(\psi)$.}
Kurzweil's theorem \cite{Kur55} {states that if $A$ is badly
approximable, then $W_A(\psi)$ has zero or full measure according as
the series \eqref{eq:classical-sum} converges or diverges.
Moreover, this law holds for every non-increasing $\psi$ only if $A$
is badly approximable.}
{Beyond badly approximable matrices, general exact criteria
were previously available only in the case $m=n=1$:}
Fuchs and Kim \cite{FK16} gave an exact
{criterion based on the continued fraction expansion; later,} Kim, Rams and
Wang \cite{KRW18} determined the exact Hausdorff dimensions of
these sets. For nonsingular matrices, the first-named author
\cite{Kim24} gave a sufficient full-measure condition, extended to
the weighted setting by Hussain and Ward \cite{HW24}.
The asymptotic problem {we consider is} to determine when $W_A(\psi)$
is null or conull for an arbitrary fixed $A$.

{We next recall the main metric results for the uniform
sets $D_A(\psi)$.}
Transference relates the approximation of shifts to homogeneous
approximation of $A^{\mathsf T}$: Jarn\'ik \cite{Jar39} gave
conditions for every shift to be uniformly approximable{,
proving an analogue of Dirichlet's theorem in the twisted setup.}
Bugeaud and Laurent \cite{BL05} related almost-everywhere
asymptotic exponents to uniform homogeneous exponents, and
Moshchevitin and the second-named author \cite{MN25} obtained
{a zero--one law for the Lebesgue measure of $D_A(\psi)$,
up to multiplicative constants in $\psi$.}
These results do not in general determine the Hausdorff
dimension of $D_A(\psi)$ when it is {Lebesgue null}. Kim and Liao \cite{KL19}
{obtained Hausdorff-dimension results} for irrational rotations when
{$\psi$} is a power function.
Aggarwal \cite{Agg26} obtained dimension bounds for singular
affine forms.
The uniform problem {we consider is} to determine dimensions for
arbitrary fixed matrices.

\subsection{The two companion papers}

The two companion papers address these questions for arbitrary
fixed matrices $A$. For non-increasing $\psi$
{satisfying the regularity condition \eqref{eq:dyadic}}, the
asymptotic paper \cite{KNA26} gives a
convergence--divergence criterion deciding whether $W_A(\psi)$
is null or conull. It also derives results on
{twisted badly approximable vectors} and Hausdorff measure and dimension.

For continuous, strictly decreasing
$\psi$ {satisfying the regularity condition \eqref{eq:lower-qm}},
the uniform paper \cite{KNU26} computes the Hausdorff
dimensions of the sets
{$\bigcap\limits_{c>0}D_A(c\psi)$ and
$\bigcup\limits_{c>0}D_A(c\psi)$, and in many cases of $D_A(\psi)$.}
For powers and near powers, it determines
$\dim_H D_A(c\psi)$ for every $c>0$, except at finitely many
exponents depending only on $m,n$.

{Different denominator vectors can give the same or nearby points
in the projected orbit $\{\langle A\bq\rangle:\bq\in\Z^n\}$ when
$A$ maps their differences close to integer vectors.
Both papers express their formulae through a homogeneous lattice count
that measures this clustering.} Section~\ref{sec:results} states the precise
hypotheses and main results, and Section~\ref{sec:geometry}
explains how this count enters the two arguments. Complete
proofs are given in the companion papers.

\medskip
\noindent\textbf{Conventions.}
All balls in $\T^m$ use the quotient supremum metric, and $\dim_H$
denotes Hausdorff dimension. Write $\mathcal L^m$ for normalized
Lebesgue measure on $\T^m$. The symbols $\ll$, $\gg$ and
$\asymp$ allow positive constants depending on $m,n$, unless
additional fixed parameters are indicated.

\medskip
\noindent\textbf{Use of artificial intelligence.}
ChatGPT Codex was used in editing, proofreading and creating figures for
this article. The authors take full responsibility for all the content of
this manuscript.

\medskip
\noindent\textbf{Acknowledgements.}
The authors are grateful to Dmitry Kleinbock for useful suggestions and
discussion.

\section{Main results}
\label{sec:results}

{First, we introduce the lattice count which enters both
asymptotic and uniform results.}

\subsection{Lattice counts and the diagonal action}

To encode the Diophantine properties of $A$, associate to it the
unimodular lattice
\begin{equation}\label{eq:lattice}
 \Lambda_A:=
 \begin{pmatrix}I_m&A\\0&I_n\end{pmatrix}\Z^d
 =\{(A\bq-\bp,\bq):
       (\bp,\bq)\in\Z^m\times\Z^n\}.
\end{equation}
{For} $r,T>0$, define
\begin{equation*}
 N_A(r,T):=
 \#\bigl(\Lambda_A\cap([-r,r]^m\times[-T,T]^n)\bigr).
\end{equation*}
Thus $N_A(r,T)$ counts integer pairs with error at most $r$ and
denominator size at most $T$. Such pairs describe differences
between nearby projected orbit points, so the count records their
clustering. It includes the zero vector and is therefore always
positive.

We use successive minima relative to the fixed unit box
$B=[-1,1]^d$. For a full-rank lattice $\Lambda$, let
$\lambda_i(\Lambda)$ be the least $\lambda>0$ such that
$\lambda B$ contains $i$ linearly independent vectors of
$\Lambda$. Consider the diagonal action
\[
 a_t:=\begin{pmatrix}
        \e^{t/m}I_m&0\\
        0&\e^{-t/n}I_n
       \end{pmatrix},\qquad t\in\R.
\]
Write $\lambda_i(t):=\lambda_i(a_t\Lambda_A)$.
For $r,T>0$, put
\[
 R_{r,T}:=(r^mT^n)^{1/d},\qquad
 t_{r,T}:=\frac{mn}{d}\log\frac{T}{r}.
\]
The map $a_{t_{r,T}}$ sends the
{rectangle $[-r,r]^m\times[-T,T]^n$} to
$R_{r,T}B$, so
$N_A(r,T)=\#(a_{t_{r,T}}\Lambda_A\cap R_{r,T}B)$.
Standard lattice-counting estimates therefore give
\begin{equation}\label{eq:minima-product}
 N_A(r,T)\asymp
 \prod_{i=1}^d
 \max\left\{1,
 \frac{R_{r,T}}{\lambda_i(t_{r,T})}\right\}.
\end{equation}
This comparison is uniform in $A,r,T$ \cite{Henk02,KNA26}.
Thus the count is determined by the successive minima along
$a_t\Lambda_A$, connecting both approximation problems with
parametric geometry of numbers \cite{SS09,DFSU24}.

\subsection{A zero--one law for asymptotic approximation}

We call $\psi$ \emph{dyadically regular} if, for some $0<\varrho<1$,
\begin{equation}\label{eq:dyadic}
 \psi(2^{k+1})\le\varrho\psi(2^k)
 \qquad\text{for all sufficiently large }k.
\end{equation}
{For instance, if} $v>0$, the power function $\psi_v(T):=T^{-v}$ satisfies this condition.
The asymptotic measure law is as follows.

\begin{thm}[{\cite[Theorem~1.1]{KNA26}}]
\label{thm:asymptotic}
Let $A\in\Mat_{m,n}$ and let $\psi:\R_+\to\R_+$ be non-increasing
and satisfy \eqref{eq:dyadic}. Then
\begin{equation}\label{eq:zero-one}
 \mathcal L^m(W_A(\psi))=
 \begin{dcases}
 0,&\displaystyle\text{if}\quad
 \sum_{k=0}^\infty
 \frac{2^{kn}\psi(2^k)^m}{N_A(\psi(2^k),2^k)}<\infty,\\[7pt]
 1,&\displaystyle\text{if}\quad
 \sum_{k=0}^\infty
 \frac{2^{kn}\psi(2^k)^m}{N_A(\psi(2^k),2^k)}=\infty.
 \end{dcases}
\end{equation}
\end{thm}

The summand in \eqref{eq:zero-one} separates volume from clustering:
\[
 \frac{2^{kn}\psi(2^k)^m}{N_A(\psi(2^k),2^k)}
 =\frac{\text{Khintchine--Groshev block volume}}
        {\text{homogeneous clustering multiplicity}}.
\]
The numerator is the $k$-th dyadic term in \eqref{eq:classical-sum},
while the denominator corrects for clustering in the projected orbit. For badly approximable matrices and irrational rotations, this
recovers the Kurzweil and Fuchs--Kim criteria, respectively,
within the class of dyadically regular functions.

\subsection{Lagrange constants and badly approximable shifts}

Define the inhomogeneous $\psi$-Lagrange constant by
\[
 \mathfrak c_{A,\psi}(\bb):=
 \liminf_{|\bq|\to\infty}
 \frac{\min_{\bp\in\Z^m}|A\bq-\bp-\bb|}{\psi(|\bq|)}.
\]
Under the hypotheses of Theorem~\ref{thm:asymptotic}, this
constant is almost surely $0$ or $+\infty$ according as the
series in \eqref{eq:zero-one} diverges or converges
\cite[Corollary~2.1]{KNA26}.

At the Dirichlet exponent {$n/m$}, define
$\operatorname{Bad}_A:=
\{\bb\in\T^m:\mathfrak c_{A,\psi_{n/m}}(\bb)>0\}$.
At time $t=nk\log2$, \eqref{eq:minima-product} gives
\[
 {\Delta_k(A):=\prod_{i=1}^d
 \min\{1,\lambda_i(a_{nk\log2}\Lambda_A)\}
 \asymp N_A(2^{-kn/m},2^k)^{-1}.}
\]
The measure of badly approximable shifts is then determined as
follows.

\begin{thm}[{\cite[Proposition~2.2]{KNA26}}]
\label{thm:report-bad-targets}
Let $A\in\Mat_{m,n}$. Then
\begin{equation*}
 \mathcal L^m(\operatorname{Bad}_A)=
 \begin{dcases}
 0,&\text{if}\quad\displaystyle\sum_k\Delta_k(A)=\infty,\\[3pt]
 1,&\text{if}\quad\displaystyle\sum_k\Delta_k(A)<\infty.
 \end{dcases}
\end{equation*}
\end{thm}

This answers the question in \cite[Remark~1.8(3)]{Kim24}
and the Lebesgue-measure case of \cite[Question~1.3]{MRS24}.

Recall that $A$ is \emph{singular} if,
for every $\varepsilon>0$, \eqref{eq:dirichlet} is solvable for
all sufficiently large $T$ with $\psi(T)=\varepsilon T^{-n/m}$;
otherwise it is \emph{nonsingular}. For a nonsingular matrix,
the count $N_A(2^{-kn/m},2^k)$ is bounded
{for infinitely many $k$}, so $\sum_k\Delta_k(A)=\infty$ and
$\operatorname{Bad}_A$ is null
\cite[Corollary~2.4]{KNA26}.

\subsection{Asymptotic Hausdorff measure and dimension}
\label{subsec:report-asymptotic-hausdorff}

To describe the size of asymptotic approximation sets more
finely, we turn to Hausdorff measure. For non-increasing,
dyadically regular $\psi$ and $0<s<m$, the count gives two
sufficient conditions
\cite[Corollary~3.1]{KNA26}:
\begin{equation*}
 \mathcal H^s(W_A(\psi))=
 \begin{dcases}
 0,&\text{if}\quad\displaystyle
 \sum_{k=0}^\infty
 \frac{2^{kn}\psi(2^k)^s}{N_A(\psi(2^k),2^k)}<\infty,\\[9pt]
 \infty,&\text{if}\quad\displaystyle
 \sum_{k=0}^\infty
 \frac{2^{kn}\psi(2^k)^s}{N_A(\psi(2^k)^{s/m},2^k)}=\infty.
 \end{dcases}
\end{equation*}
The two denominators need not agree, so these implications do
not form a single convergence--divergence criterion. The first
counts clusters at the original radius, while the second uses
the enlarged radius $\psi(2^k)^{s/m}$ required by mass
transference. When the critical series $\sum_k\Delta_k(A)$
diverges, these tests nevertheless determine the Hausdorff
dimension for every power function \cite[Corollary~3.2]{KNA26}:
\begin{equation*}
 \sum_k\Delta_k(A)=\infty
 \quad\Longrightarrow\quad
 \dim_H W_A(\psi_v)=\min\left\{m,\frac nv\right\}
 \quad(v>0).
\end{equation*}
In particular, this formula holds for every nonsingular $A$.

\subsection{Uniform Hausdorff dimension}

The dimension of $D_A(c\psi)$ may depend on $c$. To
{eliminate this dependence}, define {the endpoint sets}
\begin{equation*}
 \Sing_A(\psi):=\bigcap_{c>0}D_A(c\psi),
 \qquad U_A(\psi):=\bigcup_{c>0}D_A(c\psi).
\end{equation*}
{We refer to the elements of $\Sing_A(\psi_{n/m})$ simply as
\emph{twisted singular vectors}; see \cite[Section~2.4]{Agg26} for the corresponding
fixed-matrix problem for singular affine forms.}

We call $\psi$ \emph{lower quasimultiplicative} if there are $\tau>0$
and $R_0>1$ such that
\begin{equation}\label{eq:lower-qm}
 R^{-\tau}\psi(T)\le\psi(RT)
 \qquad\text{for all large }T\text{ and every }R\ge R_0.
\end{equation}
For continuous, strictly decreasing $\psi$ tending to zero,
we {will assume} $\psi(T)\to\infty$ as $T\to0$,
so that $\psi^{-1}$ is defined on $\R_+$.
{This assumption does not affect the set $D_A(\psi)$.}

For $L,c>0$, set
\begin{equation}\label{eq:levels}
 T_j:=\psi^{-1}(\e^{-jL}),\qquad j\ge0.
\end{equation}
and
\begin{equation}\label{eq:pressure}
 \mathcal P_{A,\psi}(L,c):=
 \liminf_{N\to\infty}\frac1{NL}
 \sum_{j=0}^{N-1}
 \log\frac{N_A(c\psi(T_j),T_{j+1})}
              {N_A(c\psi(T_{j+1}),T_{j+1})}.
\end{equation}
The count ratio governs how many smaller approximation balls
can be retained at each step; $NL=-\log\psi(T_N)$ measures
the total logarithmic shrinkage of the radii.
{Our main result regarding the uniform setup is as follows.}

\begin{thm}[{\cite[Theorem~1.2]{KNU26}}]
\label{thm:uniform}
Let $A\in\Mat_{m,n}$ and let $\psi$ be continuous, strictly
decreasing, lower quasimultiplicative and tend to zero.
\begin{enumerate}[label=\textup{(\roman*)}]
\item The limits
\begin{equation*}
 \begin{split}
 P_-(A,\psi)&:=\lim_{L\to\infty}\inf_{c>0}
                  \mathcal P_{A,\psi}(L,c),\\
 P_+(A,\psi)&:=\lim_{L\to\infty}\sup_{c>0}
                  \mathcal P_{A,\psi}(L,c)
 \end{split}
\end{equation*}
exist, and
\begin{equation*}
 \dim_H\Sing_A(\psi)
 =\inf_{c>0}\dim_H D_A(c\psi)=P_-(A,\psi),
\end{equation*}
\begin{equation*}
 \dim_H U_A(\psi)
 =\sup_{c>0}\dim_H D_A(c\psi)=P_+(A,\psi).
\end{equation*}
\item The quantity $\dim_H D_A(c\psi)$ is independent of $c>0$
if and only if $P_-(A,\psi)=P_+(A,\psi)$. Under this condition,
\begin{equation*}
 \dim_H D_A(c\psi)
 =\lim_{L\to\infty}\mathcal P_{A,\psi}(L,1)
 \qquad\text{for every }c>0.
\end{equation*}
\end{enumerate}
\end{thm}

{In particular, $\dim_H D_A(\psi)$} lies between $P_-$ and $P_+$.

Under the hypotheses of Theorem~\ref{thm:uniform}, both
$\Sing_A(\psi)$ and $U_A(\psi)$ are Borel subgroups of
$H_A:=\overline{\{\langle A\bq\rangle:\bq\in\Z^n\}}$
and each has either zero or locally full $\mathcal H^f$-measure
for every dimension function $f$ \cite[Theorem~1.9]{KNU26}.
Here $f:\R_+\to\R_+$ is continuous and nondecreasing with
$f(r)\to0$ as $r\to0$; locally full means
$\mathcal H^f(E\cap B)=\mathcal H^f(B)$ for every nonempty
relative ball $B$ in $H_A$.

\subsection{Power functions and uniform exponents}
\label{subsec:report-powers}

For power functions, the dimension of
$D_A(c\psi_\tau)$ is independent of $c>0$ outside the finite set
\begin{equation*}
 \mathcal E_{m,n}:=
 \left\{\frac ba:a,b\in\N,\ 1\le a\le m,\ 1\le b\le n\right\}.
\end{equation*}
{In this case, we can state the dimension formula in a simpler
way than in Theorem~\ref{thm:uniform}. Put}
\begingroup

\begin{equation*}
 t_\tau(T):=\tfrac{mn(1+\tau)}d\log T,\quad
 r_{A,\tau}(T):=\#\{1\le i\le d:\lambda_i(t_\tau(T))<T^{(n-m\tau)/d}\}.
\end{equation*}
\endgroup
Thus $r_{A,\tau}(T)$ counts the minima below {$T^{(n-m\tau)/d}$}
at time $t_\tau(T)$ along $a_t\Lambda_A$.

\begin{thm}[{\cite[Theorem~1.4]{KNU26}}]
\label{thm:report-near-powers}
Let $A\in\Mat_{m,n}$ and
$\tau\in(0,\infty)\setminus\mathcal E_{m,n}$. For every $c>0$,
\begin{equation}\label{eq:report-power-dimension}
 \dim_H D_A(c\psi_\tau)
 =\liminf_{T\to\infty}
 \frac{\displaystyle\int_1^T r_{A,\tau}(S)\,\frac{dS}{S}
       -\log N_A(T^{-\tau},T)}{(1+\tau)\log T}.
\end{equation}
Write $\delta_A(\tau)$ for this common value. It is continuous on each
component of $(0,\infty)\setminus\mathcal E_{m,n}$.
Moreover, if $\psi$ is continuous, strictly decreasing to zero,
and
\begin{equation*}
 \lim_{T\to\infty}\frac{\log\psi(T)^{-1}}{\log T}=\tau,
\end{equation*}
then
\begin{equation*}
 \dim_H D_A(c\psi)=\delta_A(\tau)
 \qquad\text{for every }c>0.
\end{equation*}
No lower quasimultiplicativity is required for this last assertion.
\end{thm}

For these powers, both $\Sing_A(\psi_\tau)$ and
$U_A(\psi_\tau)$ have dimension $\delta_A(\tau)$.
The near-power assertion applies, for example, to
{$\psi(T)=T^{-\tau}(\log T)^\beta$} for every fixed $\beta\in\R$.

The exclusions are necessary: for each
$\tau=b/a\in\mathcal E_{m,n}$, there is a matrix $A$ with
$\dim_H D_A(c\psi_\tau)=0$ for small $c$ and
{$\dim_H D_A(c\psi_\tau)=a$} for large $c$
\cite[Remark~1.5]{KNU26}.
{In the $m=n=1$ case, for irrational $A$}, continued fractions give {formulae
for $\dim_H D_A(\psi)$} under
the hypotheses of Theorem~\ref{thm:uniform} whenever the ratio
$\log\psi(T)^{-1}/\log T$ stays uniformly below or above $1$
\cite[Theorem~2.2]{KNU26}, extending \cite{KL19}.

These {formulae} also describe the uniform exponent, defined by
$\widehat\omega(A,\bb):=\sup\{s>0:\bb\in D_A(\psi_s)\}$,
with $\sup\varnothing=0$. Write $\delta_A(\tau^\pm)$ for
the one-sided limits of the non-increasing function $\delta_A$,
taken through exponents outside $\mathcal E_{m,n}$.
For every $\tau>0$, including {those belonging to the exceptional
set $\mathcal E_{m,n}$},
\cite[Proposition~1.6]{KNU26} gives
\begin{equation*}
 \begin{split}
 \dim_H\{\bb:\widehat\omega(A,\bb)\ge\tau\}&=\delta_A(\tau^-),\\
 \dim_H\{\bb:\widehat\omega(A,\bb)>\tau\}&=\delta_A(\tau^+).
 \end{split}
\end{equation*}
Every nonempty exact level
$\{\bb:\widehat\omega(A,\bb)=\tau\}$ has dimension
$\delta_A(\tau^-)$.

\subsection{{Twisted singular vectors}}
\label{subsec:report-escape}

For singular affine forms, a simpler bound uses only the first
minimum along $a_t\Lambda_A$. Define the
\emph{lower escape of mass} by
\begin{equation*}
 \underline{\Emass}(A):=
 \lim_{\varepsilon\searrow0}\liminf_{X\to\infty}
 \frac1X\int_0^X
 \mathbf1_{\{\lambda_1(t)<\varepsilon\}}\,dt.
\end{equation*}
It measures the lower proportion of time spent with an
arbitrarily small first minimum.

\begin{thm}[{\cite[Propositions~2.4--2.5]{KNU26}}]
\label{thm:report-escape}
For every $A\in\Mat_{m,n}$,
\begin{equation*}
 \dim_H\Sing_A(\psi_{n/m})
 \le\frac{m(d-1)}d\,\underline{\Emass}(A).
\end{equation*}
For every $\eta\in[0,1]$, there is an $A\in\Mat_{m,n}$ with
$\underline{\Emass}(A)=\eta$ attaining equality.
\end{thm}

This improves the coefficient $m$ in \cite{Agg26}; {the theorem
is optimal} when only lower escape of mass
is used.

\section{The common geometry: clusters and scales}
\label{sec:geometry}

\begin{figure}[htbp]
 \centering
\begin{tikzpicture}[
  x=1cm,y=1cm,>=Stealth,font=\footnotesize,text=figink,
  paneltitle/.style={font=\small\bfseries,align=center},
  orbit/.style={fill=figorbit,draw=white,line width=.18pt}
]
  \path[use as bounding box] (0,-1.25) rectangle (11.6,4.78);
  \node[paneltitle] at (1.95,4.55) {Partition $\mathcal P(T)$};
  \node[paneltitle] at (9.65,4.55) {Good cells $\mathcal G(T)$};
  \node at (1.95,4.23)
    {Blue points: $\langle A\mathbf q\rangle$, $|\mathbf q|\le T$};

  \foreach \xbase in {.45,8.15}{
    \fill[black!2] (\xbase,1.05) rectangle +(3,3);
  }
  \foreach \i/\j in {0/1,1/4,2/2,3/2,4/4,4/0}{
    \fill[orange!18] ({8.15+.5*\i},{1.05+.5*\j}) rectangle +(.5,.5);
  }
  \foreach \xbase in {.45,8.15}{
    \foreach \k in {1,...,5}{
      \draw[draw=figline!65,line width=.35pt]
        ({\xbase+.5*\k},1.05)--+(0,3)
        (\xbase,{1.05+.5*\k})--+(3,0);
    }
    \draw[draw=figline,line width=.65pt]
      (\xbase,1.05) rectangle +(3,3);
    \foreach \u/\v in {
      .13/.66,.28/.82,.39/.68,
      .63/2.14,.78/2.31,.88/2.19,
      1.12/1.14,1.27/1.29,1.36/1.12,1.23/1.40,
      1.62/1.18,1.76/1.34,1.89/1.23,
      2.10/2.18,2.28/2.35,2.39/2.12,
      2.11/.16,2.29/.33,2.39/.18}{
      \filldraw[orbit] ({\xbase+\u},{1.05+\v}) circle (1.35pt);
    }
  }

  \draw[->,draw=figmuted,line width=.8pt] (3.89,3.02)--(7.68,3.02);
  \node[align=center] at (5.79,1.62)
    {{Number of points in cluster:}\\$N_A(\psi(T),T)$};
  \draw[->,draw=figwindow,line width=.65pt]
    (7.20,1.30)--(8.35,1.78);

  \draw[draw=figwindow,line width=.5pt]
    (.45,.82)--(.95,.82)
    (.45,.76)--(.45,.88) (.95,.76)--(.95,.88);
  \node[below=2pt,text=figwindow] at (.70,.82) {$\asymp\psi(T)$};
  \node[align=center] at (2.40,-.05)
    {{Total number of points} $\asymp T^n$};
  \node[align=center] at (5.80,-.80)
    {$\displaystyle\#\mathcal G(T)\approx
      \frac{\text{{total number of points}}}
           {\text{{number of points in one cell}}}
      \approx\frac{T^n}{N_A(\psi(T),T)}$};
\end{tikzpicture}

 \caption{Blue orbit centers visit the orange good cells.}
 \label{fig:partition}
\end{figure}
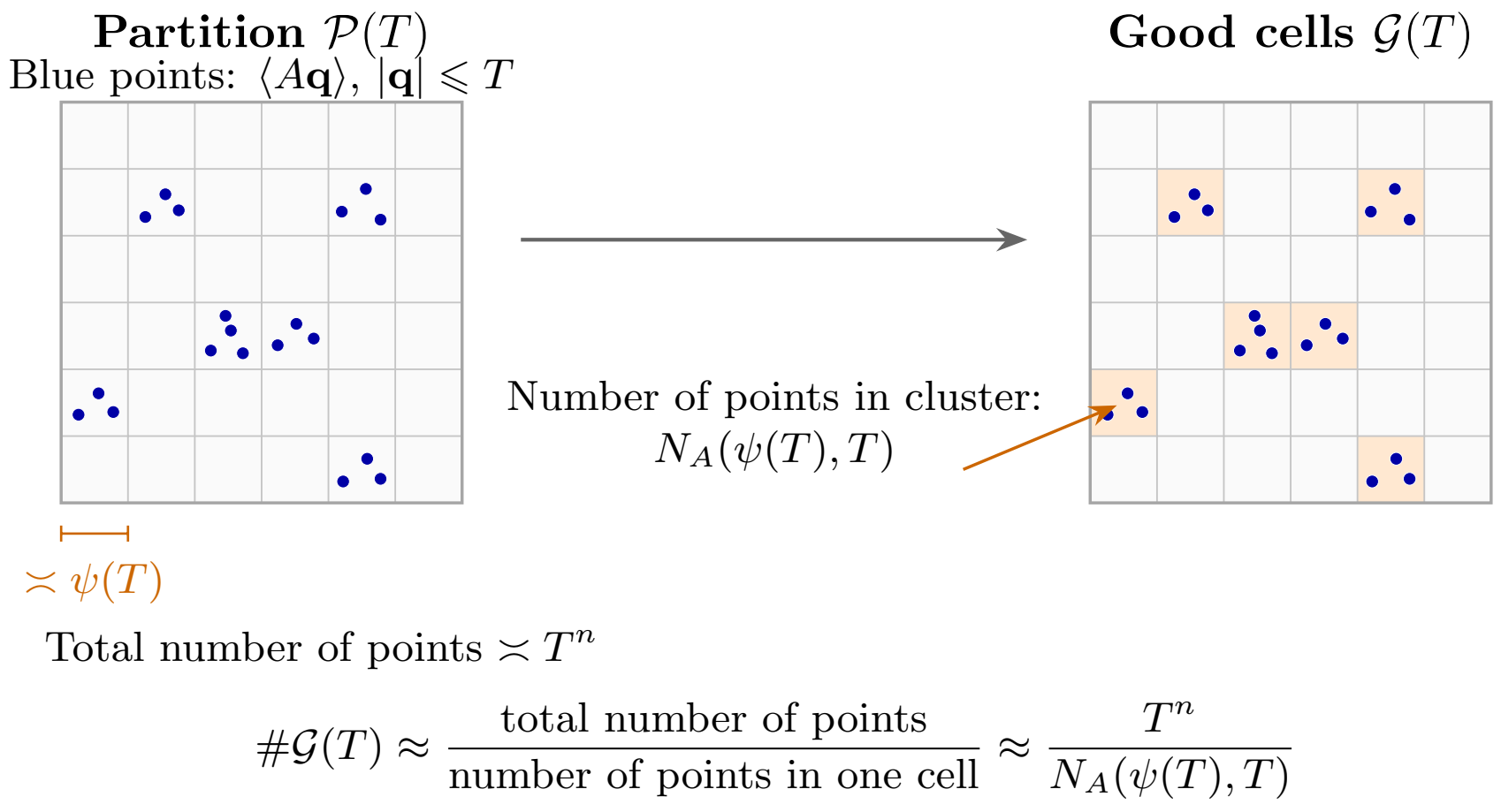

Fix a large $T$. Partition $\T^m$ into
half-open cubes of side comparable to $\psi(T)$ and at most $\psi(T)$;
in the supremum norm these are ball-shaped cells. Denote this
partition by $\mathcal P(T)$, and call a cell \emph{good} if it
contains at least one orbit point with denominator bound $T$:
\begin{equation*}
 \mathcal G(T):=\{B\in\mathcal P(T):
    \langle A\bq\rangle\in B\text{ for some }|\bq|\le T\},
 \qquad E(T):=\bigcup_{B\in\mathcal G(T)}B.
\end{equation*}
Every shift in a good cell satisfies \eqref{eq:system}, since it
is within $\psi(T)$ of the orbit point in that cell.
{Conversely, a shift satisfying \eqref{eq:system}
lies within distance $\psi(T)$ of a good cell. Thus the partition
description agrees with the approximation condition up to a
multiplicative change of radius, bounded independently of $A$
and $T$.}
Figure~\ref{fig:partition} shows the question we need to answer:
how many cells are good?

There are about $T^n$ {points $\langle A\bq\rangle$
with $|\bq|\le T$, counted with multiplicity}, but this need not give
$T^n$ different good cells. Several vectors can visit the same
cell, or even give exactly the same torus point. To understand
this multiplicity, lift the picture from $\T^m$ to the lattice
$\Lambda_A\subset\R^{m+n}$ in \eqref{eq:lattice}. The extra
coordinates keep track of the denominator vectors {$\bq$} rather than
forgetting them under projection.

First count all visits. For each $\bq\in\Z^n$ with $|\bq|\le T$,
there is exactly one $\bp\in\Z^m$ for which
$A\bq-\bp\in[0,1)^m$. Hence
\begin{equation*}
 \underbrace{\#\bigl(\Lambda_A\cap
   ([0,1)^m\times[-T,T]^n)\bigr)}_{\text{total counting}}
 =(2\lfloor T\rfloor+1)^n\asymp T^n.
\end{equation*}
This is comparable to the volume of the cylinder.

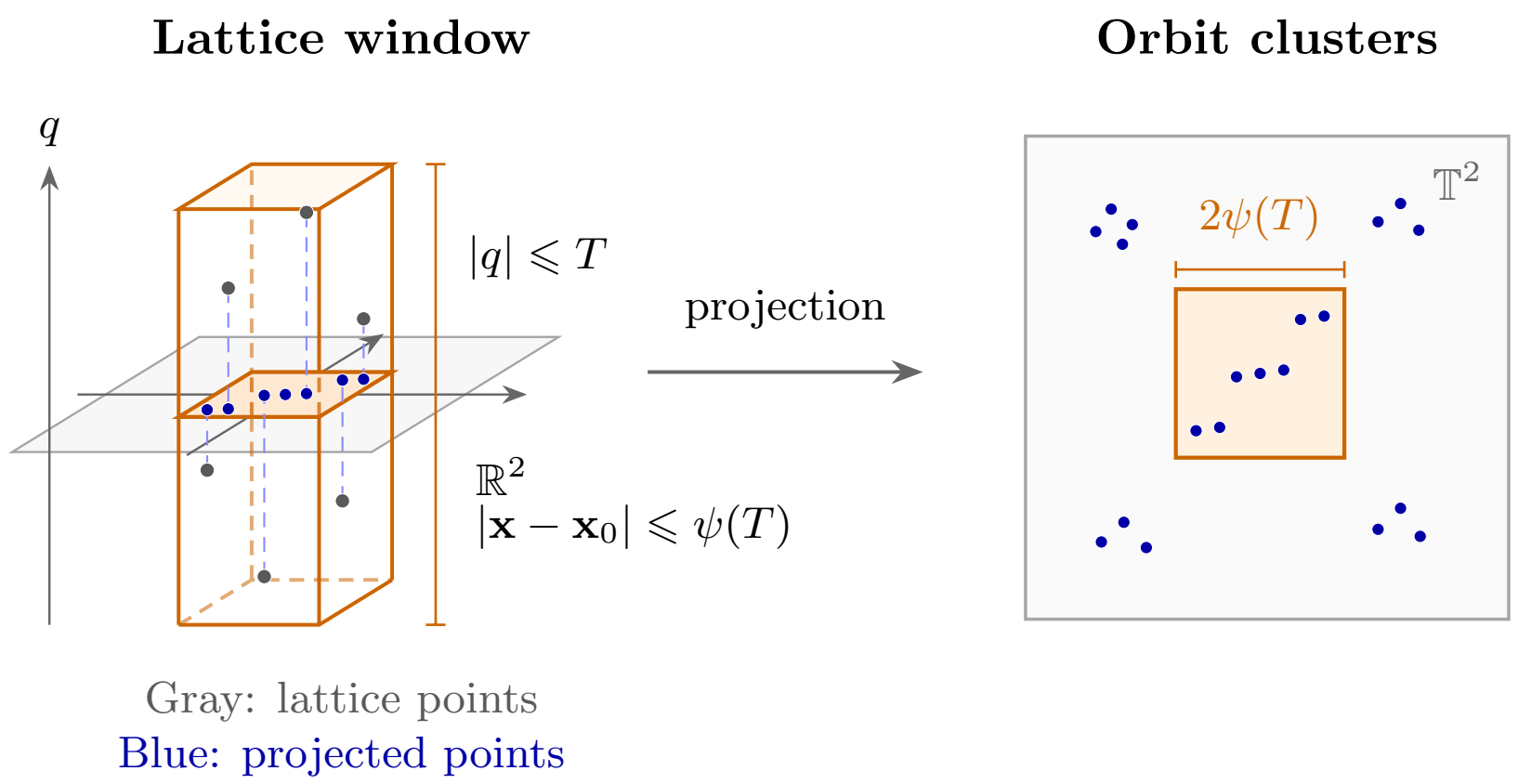
\begin{figure}[htbp]
 \centering
\begin{tikzpicture}[
  x=1cm,y=1cm,>=Stealth,font=\footnotesize,text=figink,
  paneltitle/.style={font=\small\bfseries,align=center},
  guide/.style={draw=figmuted,line width=.45pt},
  window/.style={draw=figwindow,line width=.75pt},
  latticepoint/.style={fill=black!65,draw=white,line width=.2pt},
  imagepoint/.style={fill=figorbit,draw=white,line width=.2pt}
]
  \path[use as bounding box] (0,0) rectangle (11.6,5.55);
  \node[paneltitle] at (2.55,5.29) {Lattice window};
  \node[paneltitle] at (9.14,5.29) {Orbit clusters};

  \begin{scope}[shift={(2.15,2.74)},
    x={(1cm,0cm)},y={(.52cm,.32cm)},z={(0cm,1cm)}]
    \filldraw[fill=black!3,draw=figline,line width=.45pt]
      (-1.28,-1.28,0)--(1.28,-1.28,0)
      --(1.28,1.28,0)--(-1.28,1.28,0)--cycle;
    \draw[guide,->] (-1.48,0,0)--(1.72,0,0);
    \draw[guide,->] (0,-1.35,0)--(0,1.35,0);

    \fill[orange!5] (-.5,-.5,1.48)--(.5,-.5,1.48)
      --(.5,.5,1.48)--(-.5,.5,1.48)--cycle;
    \draw[window,densely dashed,opacity=.55]
      (-.5,.5,-1.48)--(-.5,.5,1.48)
      (-.5,-.5,-1.48)--(-.5,.5,-1.48)--(.5,.5,-1.48);
    \draw[window]
      (-.5,-.5,-1.48)--(.5,-.5,-1.48)--(.5,.5,-1.48)
      (-.5,-.5,-1.48)--(-.5,-.5,1.48)
      (.5,-.5,-1.48)--(.5,-.5,1.48)
      (.5,.5,-1.48)--(.5,.5,1.48);
    \draw[window] (-.5,-.5,1.48)--(.5,-.5,1.48)
      --(.5,.5,1.48)--(-.5,.5,1.48)--cycle;

    \filldraw[fill=orange!18,draw=figwindow,line width=.8pt]
      (-.5,-.5,0)--(.5,-.5,0)--(.5,.5,0)--(-.5,.5,0)--cycle;
    \draw[window,densely dashed,opacity=.75]
      (-.5,.5,-.32)--(-.5,.5,0)
      (.5,-.5,0)--(.5,-.5,.32);

    \foreach \u/\v/\w in {.38/.34/.43,-.38/-.34/-.43,
         -.24/-.32/.86,.24/.32/-.86,.14/.02/1.29,-.14/-.02/-1.29}{
      \draw[draw=blue!45,densely dashed,line width=.4pt]
        (\u,\v,\w)--(\u,\v,0);
    }
    \foreach \u/\v/\w in {.38/.34/.43,-.38/-.34/-.43,
         -.24/-.32/.86,.24/.32/-.86,.14/.02/1.29,-.14/-.02/-1.29}{
      \filldraw[latticepoint] (\u,\v,\w) circle (1.45pt);
      \filldraw[imagepoint] (\u,\v,0) circle (1.25pt);
    }
    \filldraw[imagepoint] (0,0,0) circle (1.25pt);
  \end{scope}

  \draw[guide,->] (.47,1.1)--(.47,4.37) node[above] {$q$};
  \draw[draw=figwindow,line width=.5pt]
    (3.15,1.1)--(3.29,1.1) (3.22,1.1)--(3.22,4.38)
    (3.15,4.38)--(3.29,4.38);
  \node[anchor=west,fill=white,inner sep=2pt] at (3.37,3.72) {$|q|\le T$};
  \node[anchor=west] at (3.37,2.16) {$\mathbb R^2$};
  \node[anchor=west] at (3.37,1.80)
    {$|\mathbf x-\mathbf x_0|\le\psi(T)$};

  \draw[->,draw=figmuted,line width=.75pt] (4.73,2.9)--(6.69,2.9);
  \node[above=5pt] at (5.71,2.9) {projection};

  \filldraw[fill=black!2,draw=figline,line width=.65pt]
    (7.42,1.14) rectangle (10.86,4.58);
  \node[anchor=north east,text=figmuted,inner sep=5pt]
    at (10.86,4.58) {$\mathbb T^2$};
  \filldraw[fill=orange!12,draw=figwindow,line width=.8pt]
    (8.49,2.29) rectangle (9.69,3.49);
  \foreach \u/\v in {.38/.34,-.38/-.34,-.24/-.32,
                       .24/.32,.14/.02,-.14/-.02,0/0}{
    \filldraw[imagepoint] ({9.09+1.20*\u},{2.89+1.20*\v}) circle (1.25pt);
  }
  \foreach \u/\v in {7.92/3.90,8.03/4.06,8.18/3.95,8.11/3.81,
      9.93/3.97,10.09/4.10,10.22/3.91,
      7.96/1.69,8.12/1.83,8.28/1.65,
      9.93/1.78,10.09/1.93,10.23/1.73}{
    \filldraw[imagepoint] (\u,\v) circle (1.25pt);
  }
  \draw[draw=figwindow,line width=.5pt]
    (8.49,3.63)--(9.69,3.63)
    (8.49,3.57)--(8.49,3.69) (9.69,3.57)--(9.69,3.69);
  \node[above=2pt,text=figwindow] at (9.09,3.63) {$2\psi(T)$};

  \node[align=center] at (2.55,.36)
    {\textcolor{black!65}{Gray: lattice points}\\
     \textcolor{figorbit}{Blue: projected points}};
\end{tikzpicture}
 \caption{{The orange window has base centered at
 $\mathbf x_0$.} Lattice points in the orange window project to a cluster
 of orbit centers ($m=2$, $n=1$).}
 \label{fig:clustering}
\end{figure}

Next look at the visits to a good cell $B$. In a lift of $B$,
these {correspond to} the points of
$\Lambda_A\cap(B\times[-T,T]^n)$. Choose one such point
$(A\bq_0-\bp_0,\bq_0)$. Since $\Lambda_A$ is a lattice,
{the translation which moves this point to the origin
preserves $\Lambda_A$}. The {number of points in a good cell}
is then measured, up to fixed
{enlargements of the window $B\times[-T,T]^n$}, by $N_A(\psi(T),T)$
(Figure~\ref{fig:clustering}).

{Fixed enlargements are needed because the chosen lattice
point $\lambda_0=(A\bq_0-\bp_0,\bq_0)$ need not be the center of its
window. If $B$ has side length $\ell\asymp\psi(T)$, then
$\bigl(3B\times[-2T,2T]^n\bigr)-\lambda_0$ contains
$[-\ell,\ell]^m\times[-T,T]^n$, where $3B$ is the concentric triple
dilate of $B$. The enlarged window therefore contains at least
$N_A(\ell,T)\asymp N_A(\psi(T),T)$ lattice points.}
These changes cost only fixed
factors, since \eqref{eq:minima-product} gives
\begin{equation*}
 {N_A(c_1r,c_2T)\asymp_{c_1,c_2}N_A(r,T)
 \qquad(c_1,c_2>0).}
\end{equation*}
One short lattice direction gives strings of nearby visits;
several {short directions} give higher-dimensional clusters. This is why the count
involves all successive minima, rather than only the shortest
vector.

We can now read off the one-scale estimate as
\begin{equation*}
 \#\mathcal G(T)\asymp
 \frac{\text{{total number of points}}}
      {\text{{number of points in one cell}}}
 \asymp\frac{T^n}{N_A(\psi(T),T)}.
\end{equation*}
Each cell has volume comparable to $\psi(T)^m$, so
\begin{equation*}
 \mathcal L^m(E(T))\asymp \psi(T)^m\#\mathcal G(T)
 \asymp\frac{\psi(T)^mT^n}{N_A(\psi(T),T)}.
\end{equation*}
Taking $T=2^k$ gives exactly the summand in
Theorem~\ref{thm:asymptotic}: the usual block volume, divided by
its clustering multiplicity. {The same volume estimate
holds for the cover by enlarged cells.}
If these volumes are summable, Borel--Cantelli gives measure zero.
For divergence, we must also
understand how the good cells at different scales overlap.

\begin{figure}[tbp]
 \centering
\begin{tikzpicture}[
  x=1cm,y=1cm,>=Stealth,font=\footnotesize,text=figink,
  parent/.style={draw=figwindow,line width=.75pt,fill=orange!5},
  paneltitle/.style={font=\small\bfseries},
  orbit/.style={fill=figorbit,draw=white,line width=.18pt}
]
  \path[use as bounding box] (0,-.95) rectangle (11.6,5.8);
  \node[paneltitle] at (2.55,5.52) {{Points} in $B_S$};
  \node[paneltitle] at (9.14,5.52) {Good cells at {level} $T$};
  \foreach \xbase in {.80,7.39}{
    \path[parent] (\xbase,1.40) rectangle +(3.5,3.5);
  }
  \foreach \i/\j in {1/5,5/5,3/3,1/1,5/1}{
    \fill[blue!12] ({7.39+.5*\i},{1.40+.5*\j}) rectangle +(.5,.5);
  }
  \foreach \k in {1,...,6}{
    \draw[draw=figline!65,line width=.35pt]
      ({7.39+.5*\k},1.40)--+(0,3.5)
      (7.39,{1.40+.5*\k})--+(3.5,0);
  }
  \foreach \xbase in {.80,7.39}{
    \node[anchor=north west,inner sep=3pt] at (\xbase,4.9) {$B_S$};
    \foreach \u/\w in {
      .59/2.62,.70/2.76,.84/2.65,.77/2.57,
      2.57/2.66,2.70/2.80,2.84/2.67,2.73/2.59,
      1.61/1.62,1.75/1.76,1.88/1.65,1.81/1.56,1.96/1.80,
      .66/.58,.80/.70,.95/.59,2.56/.57,2.70/.72,2.85/.59}{
      \filldraw[orbit] ({\xbase+\u},{1.40+\w}) circle (1.4pt);
    }
  }
  \draw[->,draw=figmuted,line width=.8pt] (4.70,3.20)--(6.99,3.20);
  \node[align=center,anchor=north] at (2.55,1.15)
    {Total\\$\approx N_A(\psi(S),T)$};
  \node[align=center,anchor=north] at (9.14,1.15)
    {Local\\$\approx N_A(\psi(T),T)$};
  \node at (5.845,-.45)
    {$\displaystyle\text{{number of level-$T$ good cells}}\approx
      \frac{\text{total}}{\text{local}}
      =\frac{N_A(\psi(S),T)}{N_A(\psi(T),T)}$};
\end{tikzpicture}
 \caption{Inside $B_S$, later orbit centers are grouped at scale
 $\psi(T)$.}
 \label{fig:transition}
\end{figure}
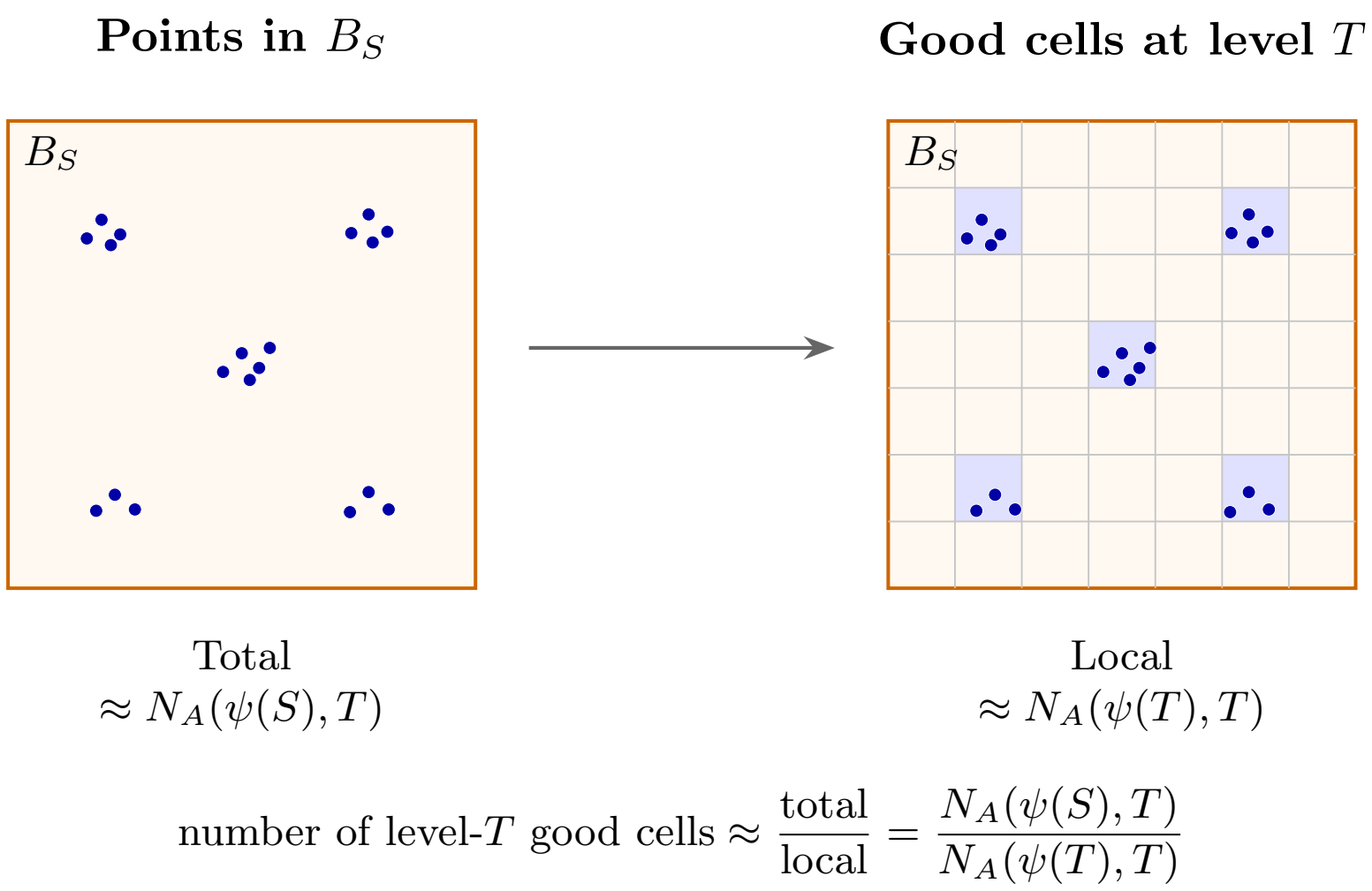

Fix $S<T$ and an earlier good cell $B_S\in\mathcal G(S)$.
At the later time $T$, the cells have {diameter comparable to} $\psi(T)$.
We ask how many
cells $B_T\in\mathcal G(T)$ meet $B_S$. The same total/local picture
now takes place inside the parent cylinder:
\[
 \underbrace{\Lambda_A\cap(B_S\times[-T,T]^n)}_{
    \text{total visits in the parent}},
 \qquad
 \underbrace{\Lambda_A\cap(B_T\times[-T,T]^n)}_{
    \text{local visits in a finer cell}}.
\]
{To count the $T$-level cells intersecting an
$S$-level cell, center the parent $B_S$ at one of its orbit points
and give it radius $\psi(S)$. For the lower bound, we also assume
$T\ge2S$ and $\psi(T)\le\psi(S)/4$.
Translation by this lattice point gives}
\begingroup

\begin{equation}\label{eq:transition-count}
 \begin{gathered}
 \text{number of level-$T$ good cells contained in }B_S\\[3pt]
 \asymp\frac{\text{number of points in one level-$S$ cell}}
              {\text{number of points in one level-$T$ cell}}
 \asymp\frac{N_A(\psi(S),T)}{N_A(\psi(T),T)}.
 \end{gathered}
\end{equation}
\endgroup
{The upper bound for cells meeting the original parent holds for
every $S<T$. Both counts use the later denominator bound $T$
(Figure~\ref{fig:transition}).}

\medskip
\noindent{\textbf{The asymptotic problem.}
The convergence case follows from the volume estimate and
Borel--Cantelli, as shown above.}
{For divergence, put $T_k=2^k$, $E_k=E(T_k)$, and}
\begingroup

\begin{equation*}
 M_k:=\frac{\psi(T_k)^mT_k^n}{N_A(\psi(T_k),T_k)},
 \qquad S_N:=\sum_{k=k_0}^N M_k,
\end{equation*}
\endgroup
{where $k_0$ is chosen so that
$\psi(T_{k+1})\le\varrho\psi(T_k)$ for $k\ge k_0$,
with $0<\varrho<1$. One has
$S_N\asymp\sum_{k=k_0}^N\mathcal L^m(E_k)$.
We seek the following bound on total pairwise overlaps, which
would show quasi-independence of the sets $E_k$ on average:}
\begingroup

\begin{equation}\label{eq:quasi-independence}
 \sum_{k,\ell=k_0}^N\mathcal L^m(E_k\cap E_\ell)
 \ll S_N^2+S_N.
\end{equation}
\endgroup
{\emergencystretch=3em
Suppose $S_N\to\infty$. A generalization of the divergent part of the
Borel--Cantelli lemma, known as the Kochen--Stone inequality
(see \cite[Lemma~5]{Spr79}), yields that $\mathcal L^m(W_A(\psi))>0$.
A separate statement proved in the asymptotic paper shows that
this measure is either zero or one; hence it equals one.\par}

{It remains to prove the quasi-independence
estimate \eqref{eq:quasi-independence}.}
Counting later good cells near each earlier good cell gives
\begingroup

\begin{equation*}
 \begin{aligned}
 \mathcal L^m(E_k\cap E_\ell)
 &\ll \underbrace{\frac{T_k^n}{N_A(\psi(T_k),T_k)}}_{
       \substack{\text{number of earlier}\\\text{good cells}}}
       \underbrace{\frac{N_A(\psi(T_k),T_\ell)}{N_A(\psi(T_\ell),T_\ell)}}_{
       \substack{\text{later good cells meeting}\\\text{one earlier cell}}}
       \underbrace{\psi(T_\ell)^m}_{
       \substack{\text{volume of one}\\\text{later cell}}}\\[3pt]
 &=M_kM_\ell\,
       \frac{N_A(\psi(T_k),T_\ell)}{\psi(T_k)^mT_\ell^n}
 \qquad(k<\ell).
 \end{aligned}
\end{equation*}
\endgroup

{
We split the pairs $k<\ell$ according to whether
$$
\Lambda_A\cap
\bigl([-\psi(T_k),\psi(T_k)]^m\times[-T_\ell,T_\ell]^n\bigr)
$$
spans $\mathbb R^{m+n}$. In the full-rank case, \eqref{eq:minima-product} shows that the count is comparable to the window volume, giving a total overlap contribution $O(S_N^2)$. In case of any smaller rank, the estimates in \cite{KNA26}, using the dyadic regularity of $\psi$, give a total contribution $O_{\varrho}(S_N)$. We add the diagonal contribution $O(S_N)$ and prove \eqref{eq:quasi-independence}.
}

\medskip
\noindent{\textbf{The uniform problem.}}
\begingroup

\emergencystretch=1em
We use the full approximation targets
\[
 F_c(T):=\bigcup_{\substack{\bq\in\Z^n\\|\bq|\le T}}
 B(\langle A\bq\rangle,c\psi(T)),\qquad c>0.
\]
Fix $c$ and choose the levels $T_j$ from \eqref{eq:levels}, so
$\psi(T_j)=\e^{-jL}$. For large $L$, lower quasimultiplicativity
ensures $T_{j+1}\ge2T_j$ for all sufficiently large $j$.

Starting at such an index $k$, \eqref{eq:transition-count} gives
an inner construction with disjoint children inside each parent
and an outer covering construction. Enlarging the balls in the
outer construction by a fixed factor makes its levels nested.
The two intersections of nested finite unions satisfy
\[
 \underline K\subseteq\bigcap_{j\ge k}F_c(T_j)
 \subseteq\overline K.
\]
In both constructions, the number of chosen children of each
level-$j$ parent is comparable to
\[
 \mathcal R_j(c):=\frac{N_A(c\psi(T_j),T_{j+1})}
                 {N_A(c\psi(T_{j+1}),T_{j+1})}.
\]
Their radii at level $N$ are comparable to $c\psi(T_N)$.
Let $C\ge1$ be a constant, independent of the level $j$ and of
sufficiently large $L$, that bounds both the ratio of the number
of chosen children to $\mathcal R_j(c)$ and its reciprocal.

\begin{itemize}[itemsep=4pt,topsep=4pt]
\item \textbf{Upper bound.}
The level-$N$ cover of $\overline K$ has at most a constant
multiple of $C^{N-k}\prod_{j=k}^{N-1}\mathcal R_j(c)$ balls.
Its $s$-dimensional covering cost is therefore at most a constant
multiple of
\[
 \underbrace{C^{N-k}\prod_{j=k}^{N-1}\mathcal R_j(c)}_{
   \substack{\text{Number of balls in the}\\\text{step-$N$ covering}}}
 \, (c\psi(T_N))^s.
\]
Taking logarithms and using $-\log\psi(T_N)=NL$ shows that
this tends to zero along a subsequence when
$s>\mathcal P_{A,\psi}(L,c)+\log C/L$. Hence
$\dim_H\overline K\le\mathcal P_{A,\psi}(L,c)+O(L^{-1})$.

\item \textbf{Lower bound.}
We use the mass distribution principle to obtain a lower bound
for $\dim_H\underline K$. Split each parent's mass equally among its
children. Each level-$N$ ball has mass at most
$C^{N-k}/\prod_{j=k}^{N-1}\mathcal R_j(c)$.
Disjointness bounds how many such balls meet a test ball of
comparable radius. The mass distribution estimate for this
construction gives
$\dim_H\underline K\ge\mathcal P_{A,\psi}(L,c)-O(L^{-1})$.
\end{itemize}
In both bounds, fixed comparison factors contribute $O(N)$ to
the logarithm, hence $O(L^{-1})$ after division by $NL$.

These estimates concern the discrete intersection
$\bigcap_{j\ge k}F_c(T_j)$. To pass to
all sufficiently large $T$, observe that
\[
 D_A(c\psi)=\bigcup_k\bigcap_{T\ge T_k}F_c(T)
 \subseteq\bigcup_k\bigcap_{j\ge k}F_c(T_j).
\]
Countable stability gives
$\dim_H D_A(c\psi)\le\mathcal P_{A,\psi}(L,c)+O(L^{-1})$.

For the lower bound, if $T_j\le T<T_{j+1}$, the denominator
chosen at $T_j$ is still at most $T$, and its error is at most
\[
 c\psi(T_j)\le\rho c\psi(T_{j+1})\le\rho c\psi(T),
 \qquad\rho=\e^L.
\]
Thus
\[
 \underline K\subseteq\bigcap_{j\ge k}F_c(T_j)
 \subseteq D_A(\rho c\psi),
\]
giving
$\dim_H D_A(\rho c\psi)\ge\mathcal P_{A,\psi}(L,c)-O(L^{-1})$.

The upper and lower bounds apply to $D_A(c\psi)$ and
$D_A(\rho c\psi)$, respectively. The discrepancy in constants
disappears for $\Sing_A(\psi)$ and $U_A(\psi)$, since their
definitions take intersections and unions over all positive
constants. Countable stability gives the dimension of $U_A(\psi)$.
For $\Sing_A(\psi)$, this gives only the upper bound
$\dim_H\Sing_A(\psi)\le P_-(A,\psi)$, since an intersection can
have smaller dimension. A separate no-dimension-drop argument
\cite[Proposition~5.4]{KNU26} proves equality, yielding the formula
for $\Sing_A(\psi)$ in Theorem~\ref{thm:uniform}.\par
\endgroup

\medskip
\noindent{\textbf{Applications.}
For asymptotic approximation, the Hausdorff-measure
tests in Section~\ref{subsec:report-asymptotic-hausdorff} follow from the same geometric considerations. The convergence test uses the $s$-dimensional covering cost of each good cell; the divergence test enlarges the targets and applies mass transference \cite{BV06}.}

{For uniform approximation, the successive-minima
expression \eqref{eq:minima-product} yields the rank formula
\eqref{eq:report-power-dimension} of
Theorem~\ref{thm:report-near-powers} for powers outside
$\mathcal E_{m,n}$. The critical power $\psi(T)=T^{-n/m}$,
for which the approximation boxes have constant volume
($\psi(T)^mT^n=1$), is treated separately: an upper estimate
using the small successive minima gives the escape-of-mass
bound in Theorem~\ref{thm:report-escape}.}

\providecommand{\bysame}{\leavevmode\hbox to3em{\hrulefill}\thinspace}

\end{document}